\documentclass[10pt,a4paper]{article}
\usepackage{amsmath}
\usepackage{amsfonts}
\usepackage{amssymb} 
\usepackage{color}
\usepackage{hyperref}
\newtheorem{theorem}{Theorem}
\newtheorem{defn}{Definition}
\newtheorem{lm}{Lemma}
\newtheorem{prop}{Proposition}

\newtheorem{rem}{Remark}
\newtheorem{plm}{Problem}

\usepackage[left=2.2cm,right=2.0cm,top=2.5cm,bottom=3cm]{geometry}
\begin{document}
	\begin{center}
		\textbf{Existence and convergence of commutative mappings on some results of fixed point theory in a class of generalized non-expansive mappings }\\\medskip
		Gezahegn Anberber Tadesse
		 Department of Mathematics, College of Computational and Natural Science, Addis Ababa University, Ethiopia.\\
	\end{center}
	\section*{Abstract} In this paper, we introduce a commutative mappings satisfying the class of generalized nonexpansive mappings which is wider than the class of mappings satisfying the condition (C), so called Condition $B_{\gamma,\mu}$. The results obtained in this paper extend and generalized non-expansive mappings and other results in this direction.
	Different properties and some fixed point results for the mappings are obtained here.\\
	\textbf{Keywords:} Nonexpansive mappings, Condition (C), Condition $B_{\gamma,\mu}$, Commuting mappings, Fixed point. 
	\section{Introduction}
	Fixed point theory and its application played an important role in many areas of applied science and solved many problems rising in engineering, mathematical economics and optimization. 
	The generalization of nonexpansive mappings and the study of related fixed point theorems with different practical applications in nonlinear functional analysis have found great importance during the recent decades \cite{afshari2016existence,aksoy2016fixed,aksoy2017fixed,alqahtani2018non,betiuk2015fixed,mishra2015fixed, dhompongsa2009edelstein, garcia2011fixed, khamsi2015monotone, gunaseelan2019coupled,mishra2015fixed, mishra2015unique, mishra2007some, patir2018fixed, sangago2010weak,  thakur2015convergence}.
	Several prominent authors \cite{agarwal2016existence,kanas2000linear,lael2016fixed,gunaseelan2019coupled,mishra2007some,patir2018fixed,suzuki2008fixed} have contributed immensely in this field, and different new classes of mappings with interesting properties have been developed in this context.
	
	Suzuki \cite{suzuki2008fixed} introduce the concept of generalized nonexpansive mappings which is called Suzuki generalized nonexpansive mappings or the condition (C) and received some fixed point results and convergence results for such mapping in Banach space. Recently in 2018,  Patir et al. \cite{patir2018some}  defined a class of generalized nonexpansive mappings on a nonempty subset C of a Banach space X which is wider than Suzuki's generalized nonexpansive mappings. Such type of mappings was called the class of mappings satisfying the condition $B_{\gamma,\mu}$ on $C$. (see \cite{patir2018some}). For a nonempty bounded and convex subset C, every self-mapping T on C satisfying $B_{\gamma,\mu}$ condition has an almost fixed point sequence \cite{garcia2011fixed}, but preserving their fixed point properties. 
	
	Having these inspiring background, we can ask:
	\begin{plm}
		Is that extend the existence and convergece theorem of the self mappings $T$ of the condition $B_{\gamma,\mu}$ to the  finite commutative mappings on this condition?
	\end{plm}
	 To answer  this, we can apply some facts related to the condition $B_{\gamma,\mu}$ on a non empty subset of a Banach space $X$.
	\section{Preliminaries}
	We now  recall the following definitions, propositions and lemmas to be use the main result.
	\begin{defn}\label{definition1}
		Let C be a nonempty subset of a Banach space $X$. A mapping
		$T : C\rightarrow  X$ is said to be \textbf{non-expansive} if $||Tx - Ty|| \leq ||x - y||$ for all $x, y \in C$.
	\end{defn}
	\begin{defn}\label{definition2}
		For a nonempty subset $C$ of a Banach space $X$, a mapping $T :C\rightarrow X$ is called \textbf{quasi-nonexpansive} if~ $||Tx-z||\leq ||x-z||$ for all $x \in C$ and $z \in F(T)$ (where
		$F(T)$ denotes the set of all fixed points of $T$).
	\end{defn}
	\begin{defn}\label{definition3}
		( \cite{suzuki2008fixed}) For a nonempty subset C of a Banach space X, a mapping $T :C\rightarrow X$ is said to satisfy \textbf{the condition $(C)$} on $C$ if~ $ \frac{1}{2}||x-Tx||\leq||x-y|| $ implies $||Tx-Ty||\leq ||x - y||$ for all $x, y \in C$.
	\end{defn}
	\begin{rem}\label{remark1}
		Every nonexpansive mapping satisfies the condition (C) on C. But there are also some noncontinuous mappings satisfying the condition $C$,   (see \cite{suzuki2008fixed}).	
	\end{rem}	
	\begin{defn}\label{definition4}
		( \cite{suzuki2008fixed}) For a nonempty subset C of a Banach space X and $\lambda \in (0, 1$), a mapping $T :C\rightarrow X$ is said to satisfy \textbf{$(C_{\lambda}$)-condition} on C if $\lambda||x - Tx|| \leq ||x - y||$ implies $||Tx - Ty|| \leq ||x - y||$ for all $x, y \in C$.
	\end{defn}
	\begin{defn}\label{definition5}
		(\cite{patir2018some}) Let C be a nonempty subset of a Banach space X. Let 
		$\gamma \in [0, 1]$ and 
		$\mu \in[0, \frac{1}{2}]$ such that $2\mu \leq \gamma$ . 
		A mapping $T : C\rightarrow X$ is said to satisfy \textbf{the condition $B_{\gamma,\mu}$ }on C if, for all $x, y$ in C, $\gamma ||x - Tx|| \leq ||x - y|| + \mu||y - Ty||  $ implies $ ||Tx - Ty|| \leq(1 - \gamma )||x - y|| + \mu(||x - Ty|| + ||y - Tx||). $
	\end{defn}
	Since, this class includes the class of nonexpansive mappings (for $ \gamma=\mu=0 $).
	\begin{rem}\label{remark2}
		Let $C$ be a nonempty subset of a Banach space $X$. If a mapping satisfies the condition (C), then it will satisfy the condition $B_{\gamma,\mu}$  for $ \gamma=\mu=0 $.
	\end{rem}
  \begin{defn}
  	
  	Let $C$ be a nonempty closed convex and bounded subset of a Banach space $X$ and a self mapping $T$ on $C$ is satisfies condition $B_{\gamma,\mu}$ on $C$, then there exists a sequence $\{x_{n}\}$ such that  $\lim_{n\rightarrow\infty}||Tx_n-x_{n}||=0$. Such a sequence is called almost ﬁxed point sequence for $T$.
  \end{defn}
	In 2005, Suzuki (\cite{suzuki2005strong}) proved the mappings $T_{1}$ and $T_{2}$ are non-expansive mapping with $T_{1}\circ T_{2}=T_{2}\circ T_{1}$ by \textbf{Lemma 3.1} have a common fixed point of $T_{1}$ and $T_{2}$. In similar way to extend the prove of finite commutative mappings for some class of generalized non-expansive have a common fixed point.

	The following are some basic properties of mappings which satisfy the condition $ B_{\gamma,\mu} $ on C.	
	\begin{prop}\label{proposition1}
		(\cite{patir2018some}) Let C be a nonempty subset of a Banach space X. Let $T : C\rightarrow C$ satisfy the condition $ B_{\gamma,\mu} $ on C. Then, for all $x, y \in C$ and for $\theta\in  [0, 1]$,
		\begin{enumerate}
			\item[(i)] $ ||Tx-T^{2}x||\leq||x-Tx||, $
			\item[(ii)] at least one of the following ((a) and (b)) holds:
			\begin{enumerate}
				\item[(a)] $ \frac{\theta}{2}||x-Tx||\leq||x-y||$,
				\item[(b)]  $ \frac{\theta}{2}||Tx-T^{2}x||\leq||Tx-y||$.\\
				The condition (a) implies $||Tx - Ty|| \leq(1 - \frac{\theta}{2}
				)||x - y|| + \mu(||x - Ty|| + ||y - Tx||)$ and the condition (b) implies $||T^{2}x - Ty|| \leq(1 - \frac{\theta}{2}
				)||Tx - y|| + \mu(||Tx - Ty|| + ||y - T^{2}x||)$.
			\end{enumerate}
			\item[(iii)] 
			\begin{align*}||x - Ty|| \leq(3-\theta)||x-Tx|| +(1-\frac{\theta}{2})||x-y||+\mu(2||x-Tx||\\+||x-Ty||+||y-Tx|| +2||Tx - T^{2}x||). 
			\end{align*}
		\end{enumerate}
	\end{prop}
	\begin{prop}\label{proposition2}
		(\cite{patir2018some}) Let C be a nonempty convex and bounded subset of a Banach space X and T be a self-mapping  on C. We assume that T satisfies the condition $ B_{\gamma,\mu} $ on C. For $x_{0} \in C$, let a sequence $\{x_{n}\}$ in C be defined by;
		\begin{equation}\label{equation1}
			x_{n+1}=\lambda Tx_{n}+(1-\lambda)x_{n},
		\end{equation}
		where $\lambda \in [\gamma , 1) - \{0\}, n \in \mathbb{N} \cup \{0\}$. Then $||Tx_{n} - x_{n}||\rightarrow0$ as $n \rightarrow\infty$.
	\end{prop}
	\section{Main results}
	In this section, we prove our main results. That is, we prove the convergence and the existence theorem of fixed point of commuting mappings satisfying the condition $B_{\gamma,\mu}$.
 First, we give the following lemma, which will play an important role in the sequel.
	\begin{lm}\label{lemma4}
		Let $C$ be a nonempty closed   convex  subset of a Banach space $X$. Let $T_{1}$ and $T_{2}$ be mappings with $T_{1}\circ T_{2}=T_{2}\circ T_{1}$ satisfying the condition $B_{\gamma,\mu}$ on $C$. Let $\{\alpha_{n}\}$ be a sequence in $(0,\frac{1}{2})$ converging to $0$, and let $\{x_{n}\}$ be a sequence in $C$ such that $\{x_{n}\}$  converges strongly to some $z\in C$.
	 Then $z$ is  a common fixed point of $T_{1}$ and $T_{2}$.  
	\end{lm}
	\textbf{Proof.} Since a  sequence $\{x_{n}\}$ in $C$  converges strongly to some $z\in C$.  \\
	And by the condition $B_{\gamma,\mu}$ we get
	\begin{equation}\label{equation2}
		||T_{1}x_{n}-T_{1}z||\leq(1-\gamma)||x_{n}-z||+\mu\Big(||x_{n}-T_{1}z||+||z-T_{1}x_{n}||\Big).
	\end{equation}
	Again by using (\ref{equation2}) 
	\begin{eqnarray*}
		||x_{n}-T_{1}z||&\leq&||x_{n}-T_{1}x_{n}||+||T_{1}x_{n}-T_{1}z||\\&\leq&||x_{n}-T_{1}x_{n}||+(1-\gamma)||x_{n}-z||+\mu\Big(||x_{n}-T_{1}z||+||z-T_{1}x_{n}||\Big)\\&\leq&||x_{n}-T_{1}x_{n}||+(1-\gamma)||x_{n}-z||+\mu\Big(||x_{n}-T_{1}z||+||z-x_{n}||+||x_{n}-T_{1}x_{n}||\Big).
	\end{eqnarray*} 
It implies that 
\[(1-\mu)||x_{n}-T_{1}z||\leq(1+\mu)||T_{1}x_{n}-x_{n}||+(1-\gamma+\mu)||x_{n}-z||.\] 
	So taking ${n}\rightarrow \infty$ and using Proposition \ref{proposition2}, we get \[(1-\mu)||z-T_{1}z||\leq 0,\]
	since $\mu\neq 1$. It implies that
	\begin{equation}\label{equation3}
	T_{1}z=z	
		\end{equation}
	  This shows that $z$ is a fixed point for $T_{1}$. We note that \begin{equation}\label{equation4}
		T_{1}oT_{2}z=T_{2}oT_{1}z=T_{2}z.
	\end{equation}
	We assume that $z$ is not a fixed point of $T_{2}$.  Put \[\epsilon=\dfrac{||T_{2}z-z||}{3}>0.\]
	Then there exists $m\in\mathbb{N}$ such that 
	\[||x_{m}-z||<\epsilon,~\frac{||(1-\alpha_{m})T_{1}x_{m}+\alpha_{m}T_{2}x_{m}-x_{m}||}{\alpha_{m}}< \epsilon.\]
	Since, 
	\begin{eqnarray*}
		3\epsilon=||T_{2}z-z|| &\leq&||T_{2}z-x_{m}||+||x_{m}-z||\\
		&<&||T_{2}z-x_{m}||+\epsilon\\ \Rightarrow 2\epsilon &<&||T_{2}z-x_{m}||.
	\end{eqnarray*}
	So, by using (\ref{equation4}) and Proposition\ref{proposition1},(i)  for $\gamma=\dfrac{\theta}{2}, \theta\in[0,1]$, we obtain, we obtain
	\begin{eqnarray*}
		||T_{2}z-x_{m}||&\leq & ||T_{2}z-(1-\alpha_{m})T_{1}x_{m}-\alpha_{m}T_{2}x_{m}||+||(1-\alpha_{m})T_{1}x_{m}+\alpha_{m}T_{2}x_{m}-x_{m}||\\ &\leq&
		(1-\alpha_{m})||T_{2}z-T_{1}x_{m}||+\alpha_{m}||T_{2}z-T_{2}x_{m}||+||(1-\alpha_{m})T_{1}x_{m}+\alpha_{m}T_{2}x_{m}-x_{m}||\\&=&(1-\alpha_{m})||T_{2}\circ T_{1} z-T_{1}x_{m}||+\alpha_{m}||T_{2}z-T_{2}x_{m}||\\&&+||(1-\alpha_{m})T_{1}x_{m}+\alpha_{m}T_{2}x_{m}-x_{m}||
		\\&\leq&(1-\alpha_{m})\Big[(1-\gamma)||T_{2}z-x_{m}||+\mu(||T_{2}z-T_{1}x_{m}||+||x_{m}-T_{2}z||)\Big]+\\&&\alpha_{m}\Big[(1-\gamma)||z-x_{m}||+\mu(||z-T_{2}x_{m}||+||x_{m}-T_{2}z||)\Big]+||(1-\alpha_{m})T_{1}x_{m}+\alpha_{m}T_{2}x_{m}-x_{m}||
			\\&\leq&(1-\alpha_{m})\Big[(1-\gamma)||T_{2}z-x_{m}||+\mu(||T_{2}z-T_{1}x_{m}||+||x_{m}-T_{2}z||)\Big]+\\&&\alpha_{m}\Big[(1-\gamma)||z-x_{m}||+\mu||z-x_{m}||+\mu||x_{m}-z||\Big]+||(1-\alpha_{m})T_{1}x_{m}+\alpha_{m}T_{2}x_{m}-x_{m}||
		\\&\leq&
		(1-\alpha_{m})\Big[\dfrac{1-\gamma+\mu}{1-\mu}\Big]||T_{2} z-x_{m}||+\alpha_{m}||z-x_{m}||+\\&&||(1-\alpha_{m})T_{1}x_{m}+\alpha_{m}T_{2}x_{m}-x_{m}||\\&\leq&
		(1-\alpha_{m})||T_{2} z-x_{m}||+\alpha_{m}||z-x_{m}||+||(1-\alpha_{m})T_{1}x_{m}+\alpha_{m}T_{2}x_{m}-x_{m}||	\\
		&<&(1-\alpha_{m})||T_{2}z-x_{m}||+\alpha_{m}\epsilon+\alpha_{m}\epsilon\\&=&(1-\alpha_{m})||T_{2}z-x_{m}||+2\alpha_{m}\epsilon\\&<&(1-\alpha_{m})||T_{2}z-x_{m}||+\alpha_{m}||T_{2}x_{m}-x_{m}||\\&=&||T_{2}z-x_{m}||,
	\end{eqnarray*}
	since $~\dfrac{1-\gamma+\mu}{1-\mu}\leq 1$,  which is a contradiction. Hence,  $z$ is a common fixed point of $T_{1}$ and $T_{2}~~~~~~~~~~~~~~~~~~~~~~~~~~~~~~~\blacksquare$.
	\begin{lm} \label{lm4}
		Let $C$ be a nonempty closed convex and bounded  subset of a  Banach space $X$. Let $T_{1}$ and $T_{2}$ be mappings with $T_{1}\circ T_{2}=T_{2}\circ T_{1}$ satisfying the condition $B_{\gamma,\mu}$ on $C$. Let $\{\alpha_{n}\}$ be a sequence in $[0,\dfrac{1}{2}]$ converging to $0$. 
		Define a sequence $\{x_{n}\}$ in $C$ by $x_{0}\in C$ and
		\begin{equation}\label{equation5}
			x_{n+1}=\lambda(1-\alpha_{n})T_{1}x_{n}+\lambda\alpha_{n}T_{2}x_{n}+(1-\lambda)x_{n},
		\end{equation}
		for $\lambda\in[\gamma,1)-\{0\},~n=0,1,2,\cdots$.	 Then
		\begin{equation}\label{eqn5}
		\lim_{n\rightarrow\infty}||(1-\alpha_{n})T_{1}x_{n}+\alpha_{n}T_{2}x_{n}-x_{n}|| =0.
		\end{equation}
		 	
	\end{lm}
\textbf{Proof.} Since a sequence $\{x_{n}\}$ in $C$ defined as (\ref{equation5}) and assume that $z$ is a common fixed point of $T_{1}$ and $T_{2}$,  then, by  condition $B_{\gamma,\mu}$ we get
\begin{eqnarray*}
||(1-\alpha_{n})T_{1}x_{n}+\alpha_{n}T_{2}x_{n}-x_{n}||&\leq&\lambda||(1-\alpha_{n})T_{1}x_{n}+\alpha_{n}T_{2}x_{n}-x_{n}||\\&=&||x_{n+1}-x_{n}||\\&\leq&||x_{n+1}-x_{n}||+\mu||x_{n}-T_{1}x_{n}||.	
\end{eqnarray*}
Implies that
\begin{eqnarray*}
	&&||(1-\alpha_{n+1})T_{1}x_{n+1}+\alpha_{n+1}T_{2}x_{n+1}-T_{1}x_{n}||\leq(1-\alpha_{n+1})||T_{1}x_{n+1}-T_{1}x_{n}||+\alpha_{n+1}||T_{2}x_{n+1}-T_{1}x_{n}||\\&\leq&(1-\alpha_{n+1})\Big[||T_{1}x_{n+1}-T_{1}z||+||T_{1}z-T_{1}x_{n}||\Big]+\alpha_{n+1}\Big[||T_{2}x_{n+1}-T_{2}z||+||T_{2}z-T_{1}x_{n}||\Big]\\&\leq&(1-\alpha_{n+1})\Big[(1-\gamma)||x_{n+1}-z||+\mu(||x_{n+1}-T_{1}z||+||z-T_{1}x_{n+1}||)+(1-\gamma)||x_{n}-z||+\mu(||x_{n}-T_{1}z||+||z-T_{1}x_{n}||)\Big]\\&&+~
	\alpha_{n+1}\Big[(1-\gamma)||x_{n+1}-z||+\mu(||x_{n+1}-T_{2}z||+||z-T_{2}x_{n+1}||)+(1-\gamma)||x_{n}-z||+\mu(||x_{n}-T_{2}z||+||z-T_{2}x_{n}||)\Big]
	\\&\leq&(1-\alpha_{n+1})\Big[(1-\gamma)||x_{n+1}-z||+\mu(||x_{n+1}-z||+||z-x_{n+1}||+(1-\gamma)||x_{n}-z||+\mu(||x_{n}-z||+||z-x_{n}||)\Big]+\\&&
	\alpha_{n+1}\Big[(1-\gamma)||x_{n+1}-z||+\mu(||x_{n+1}-z||+||z-x_{n+1}||)+(1-\gamma)||x_{n}-z||+\mu(||x_{n}-z||+||z-x_{n}||)\Big]
	\\&\leq&(1-\alpha_{n+1})\Big[||x_{n+1}-z||+||x_{n}-z||\Big]+
	\alpha_{n+1}\Big[||x_{n+1}-z||+||x_{n}-z||\Big]
	\\&=&
	||x_{n+1}-z||+||x_{n}-z||.	
\end{eqnarray*}
	Now from the following we have that, 
\[||(1-\alpha_{n})T_{1}x_{n}+\alpha_{n}T_{2}x_{n}-x_{n}||\leq(1-\alpha_{n})||T_{1}x_{n}-x_{n}||+\alpha_{n}||T_{2}x_{n}-x_{n}||.\]
Taking $n\rightarrow\infty$ on both side of the above inequality and using Proposition \ref{proposition2} we get
\[\lim_{n\rightarrow\infty}||(1-\alpha_{n})T_{1}x_{n}+\alpha_{n}T_{2}x_{n}-x_{n}||=0.\]
For a nonempty compact convex subset $C$ of $X$, we have the following fixed point result.
	\begin{theorem}\label{theorem1}
	Let $C$ be a compact and convex subset of a  Banach space $X$. Let $T_{1}$ and $T_{2}$ be mappings with $T_{1}\circ T_{2}=T_{2}\circ T_{1}$ satisfying the condition $B_{\gamma,\mu}$ on $C$.  Let $\{\alpha_{n}\}$ be a sequence in $[0,\frac{1}{2}]$ converging to $0$. 
	For $x_{0}\in C$, let $\{x_{n}\}$ be a sequence in $C$ as defined by (\ref{equation5}).	 Then $\{ x_{n}  \}$ converges strongly to a common fixed point $T_{1}$ and $T_{2}$.	
	\end{theorem}
\textbf{Proof.}	
Since $C$ is compact, there exists a subsequence
$\{x_{n_{j}}\}$ of $\{x_{n}\}$ and $z\in C$ such that $\{x_{n_{j}}\}$ converges to $z\in C$ By (Lemma \ref{lemma4})\\
Now by Proposition \ref{proposition1} (ii), for $\gamma=\frac{\theta}{2},~\theta\in[0,1]$,
	\begin{eqnarray*}
	&&\gamma||x_{n_{j}}-T_{1}x_{n_{j}}||\leq||x_{n_{j}}-z||\\
	&&\Rightarrow\gamma||x_{n_{j}}-T_{1}x_{n_{j}}||\leq||x_{n_{j}}-z||+\mu||z-T_{1}z||.
\end{eqnarray*}
So by condition $B_{\gamma,\mu}$
\begin{equation}\label{eqn7}
	||T_{1}x_{n_{j}}-T_{1}z||\leq(1-\gamma)||x_{n_{j}}-z||+\mu(||x_{n_{j}}-T_{1}z||+||z-T_{1}x_{n_{j}}||)
\end{equation}
Again, by using (\ref{eqn7}) we have
	\begin{eqnarray*}
	||x_{n_{j}}-T_{1}z||&\leq&||x_{n_{j}}-T_{1}x_{n_{j}}||+||T_{1}x_{n_{j}}-T_{1}z||\\&\leq&||x_{n_{j}}-T_{1}x_{n_{j}}||+(1-\gamma)||x_{n_{j}}-z||+\mu(||x_{n_{j}}-T_{1}z||+||z-T_{1}x_{n_{j}}||)
	\\&\leq&||x_{n_{j}}-T_{1}x_{n_{j}}||+(1-\gamma)||x_{n_{j}}-z||+\mu(||x_{n_{j}}-T_{1}z||+||z-x_{n_{j}}||+||x_{n_{j}}-T_{1}x_{n_{j}}||).	
\end{eqnarray*}
So taking $n_{j}\rightarrow\infty$ and using Proposition \ref{proposition2}, we get
	\begin{eqnarray*}
	&&(1-\mu)||z-T_{1}z||\leq 0\\
	&&\Rightarrow Tz=z,	
\end{eqnarray*}
since $\mu\neq1$, shows that $z$ is a fixed point $T_{1}$. Similarly by applying Lemma \ref{lemma4}, we have that $z=T_{2}z$. Hence, $z$ is a common fixed point of $T_{1}$ and $T_{2}$.\\ 
Now by the condition $B_{\gamma,\mu}$ we have 
	\begin{eqnarray*}
		||T_{1}x_{n}-z||&=&||T_{1}x_{n}-T_{1}z||\\&\leq&(1-\gamma)||x_{n}-z||+\mu(||x_{n}-T_{1}z||+||z-T_{1}x_{n})||\\&=&(1-\gamma)||x_{n}-z||+\mu||x_{n}-z||+\mu||z-T_{1}x_{n}||.
	\end{eqnarray*}
Which implies that
\begin{equation}\label{equation7}
||T_{1}x_{n}-z||\leq(\dfrac{1-\gamma+\mu}{1-\mu})||x_{n}-z||\leq||x_{n}-z||,	
\end{equation}
since $2\mu\leq\gamma$. Similarly we have that $||T_{2}x_{n}-z||\leq||x_{n}-z||$.
Now  by using (\ref{equation5}) and (\ref{equation7}) , we have that:
	\begin{eqnarray*}
		||x_{n+1}-z||&=&||\lambda(1-\alpha_{n}) T_{1}x_{n}+\lambda\alpha_{n} T_{2}x_{n}+(1-\lambda)x_{n}-z||\\&\leq&
		\lambda(1-\alpha_{n})||T_{1}x_{n}-T_{1}z||+\lambda\alpha_{n}||T_{2}x_{n}-T_{2}z||+(1-\lambda)||x_{n}-z||
		\\&\leq&
		\lambda(1-\alpha_{n})\Big[(1-\gamma)||x_{n}-z||+\mu(||x_{n}-T_{1}z||+||z-T_{1}x_{n}||)\Big]+\\&&\lambda\alpha_{n}\Big[(1-\gamma)||x_{n}-z||+\mu(||x_{n}-T_{2}z||+||z-T_{2}x_{n}||)\Big]+(1-\lambda)||x_{n}-z||
		\\&\leq&	\lambda(1-\alpha_{n})\Big[(1-\gamma)||x_{n}-z||+2\mu||x_{n}-z||\Big]+\\&&\lambda\alpha_{n}\Big[(1-\gamma)||x_{n}-z||+2\mu||x_{n}-z||\Big]+(1-\lambda)||x_{n}-z||
		\\&\leq&\lambda(1-\alpha_{n})\Big[(1-\gamma)||x_{n}-z||+\gamma||x_{n}-z||\Big]+\\&&\lambda\alpha_{n}\Big[(1-\gamma)||x_{n}-z||+\gamma||x_{n}-z||\Big]+(1-\lambda)||x_{n}-z||
		\\&=&
		\lambda(1-\alpha_{n})||x_{n}-z||+\lambda\alpha_{n}||x_{n}-z||+(1-\lambda)||x_{n}-z||\\&=&||x_{n}-z||.
	\end{eqnarray*}
	Thus, $\{||x_{n}-z||\}$ is a monotonically decreasing non-negative sequence of real numbers, and will converges to some real numbers, say $u$. \\
	Now   by using (\ref{equation6}) we get,
\begin{eqnarray*}
	||x_{n}-z||= ||x_{n}-T_{1}z||&\leq& ||x_{n}-T_{1}x_{n}||+(1-\gamma)||x_{n}-z||+\mu\Big(||x_{n}-T_{1}z||+||z-T_{1}x_{n}||\Big)\\&\leq&||x_{n}-T_{1}x_{n}||+(1-\gamma)||x_{n}-z||+\mu\Big(2||x_{n}-z||+||x_{n}-T_{1}x_{n}||\Big),
\end{eqnarray*}
since $z=T_{1}oT_{2}z=T_{2}oT_{1}z=T_{1}z=T_{2}z,$. Then it implies that
\begin{equation}\label{equation8}
	||x_{n}-z||\leq||x_{n}-T_{1}x_{n}||+(1-\gamma)||x_{n}-z||+\mu\Big(2||x_{n}-z||+||x_{n}-T_{1}x_{n}||\Big)
\end{equation} 
	Taking the limit as $n\rightarrow\infty$ of (\ref{equation8}) and using Proposition \ref{proposition2}  we get $u\leq(1-\gamma)u+\mu(2u)$ implies that $(\gamma-2\mu)u\leq0$, which is  possible only for $u=0$, since $2\mu\leq\gamma$. Hence, $\{x_{n}\}$ converges strongly to the common fixed point of $T_{1}$ and $T_{2}.~~~~~~~~~~~~~~~~~~~~~~~~~~~~~~~\blacksquare$
	\begin{lm}\label{lemma5}
		Let $C$ be a nonempty closed convex  subset of a  Banach space $X$. Let $T_{1}, T_{2}$ and $T_{3}$ be commuting  mappings   satisfying the condition $B_{\gamma,\mu}$ on $C$. Let $\{\alpha_{n}\}$ be a sequence in $(0,\dfrac{1}{2})$ converging to $0$, and let $\{x_{n}\}$ be a sequence in $C$ such that $\{x_{n}\}$  converges strongly to some $z\in C$.
	 Then $z$   a common fixed point of $T_{1}, T_{2}$ and $T_{3}$.  
	\end{lm}
	\textbf{Proof.} Assume that a  sequence $\{x_{n}\}$ in $C$  converges strongly to some $z\in C$. Since by Lemma \ref{lemma4}, we have that $z$ is a common fixed point of $T_{1}$ and $T_{2}$. And we note that;
	\begin{equation}\label{equation9}
	T_{1}\circ T_{3}z=T_{3}\circ T_{1}z=T_{3}z,~~T_{2}\circ T_{3} z=T_{3}\circ T_{2}z=T_{3}z.
	\end{equation} 
	We assume that $z$ is not a fixed point of $T_{3}$.  Put \[\epsilon=\dfrac{||T_{3}z-z||}{3}>0.\]
	Then there exists $m\in\mathbb{N}$ such that 
	\[||x_{m}-z||<\epsilon,~\frac{||(1-\alpha_{m}-\alpha^{2}_{m})T_{1}x_{m}+\alpha_{m}T_{2}x_{m}+\alpha^{2}_{m}T_{2}x_{m}-x_{m}||}{\alpha^{2}_{m}}< \epsilon,\]
	since, 
	\begin{eqnarray*}
		3\epsilon=||T_{3}z-z|| &\leq&||T_{3}z-x_{m}||+||x_{m}-z||\\
		&<&||T_{3}z-x_{m}||+\epsilon\\ \Rightarrow 2\epsilon &<&||T_{3}z-x_{m}||.
	\end{eqnarray*}
	So, by using   Proposition\ref{proposition1} for $\gamma=\dfrac{\theta}{2}, \theta\in[0,1]$ and (\ref{equation9}), we obtain
	\begin{eqnarray*}
		||T_{3}z-x_{m}||&\leq & ||T_{3}z-(1-\alpha_{m}-\alpha^{2}_{m})T_{1}x_{m}-\alpha_{m}T_{2}x_{m}-\alpha^{2}_{m}T_{3}x_{m}||\\
		&&+||(1-\alpha_{m}-\alpha^{2}_{m})T_{1}x_{m}+\alpha_{m}T_{2}x_{m}+\alpha^{2}_{m}T_{3}x_{m}-x_{m}||\\
		&\leq&(1-\alpha_{m}-\alpha^{2}_{m})||T_{3}z-T_{1}x_{m}||+\alpha_{m}||T_{3}z-T_{2}x_{m}||+\alpha^{2}_{m}||T_{3}z-T_{3}x_{m}||\\
		&&+||(1-\alpha_{m}+\alpha^{2}_{m})T_{1}x_{m}+\alpha_{m}T_{2}x_{m}-\alpha^{2}_{m}T_{3}x_{m}-x_{m}||\\&=&(1-\alpha_{m}-\alpha^{2}_{m})||T_{3}\circ T_{1}z-T_{1}x_{m}||+\alpha_{m}||T_{3}\circ T_{2}z-T_{2}x_{m}||+\alpha^{2}_{m}||T_{3}z-T_{3}x_{m}||\\
		&&+||(1-\alpha_{m}+\alpha^{2}_{m})T_{1}x_{m}+\alpha_{m}T_{2}x_{m}-\alpha^{2}_{m}T_{3}x_{m}-x_{m}||\\
		&\leq&
		(1-\alpha_{m}-\alpha^{2}_{m})\Big[(1-\gamma)||T_{3}z-x_{m}||+\mu(||T_{3}x-x_{m})||+||x_{m}-T_{3}z||\Big]\\&&+\alpha_{m}\Big[(1-\gamma)||T_{3}z-x_{m}||+\mu(||T_{3}z-T_{2}x_{m}||+||x_{m}-T_{3}z||)\Big]+\alpha^{2}_{m}||T_{3}z-T_{3}x_{m}||\\
		&&+||(1-\alpha_{m}+\alpha^{2}_{m})T_{1}x_{m}+\alpha_{m}T_{2}x_{m}-\alpha^{2}_{m}T_{3}x_{m}-x_{m}||\\
		&\leq&
		(1-\alpha_{m}-\alpha^{2}_{m})\Big[\dfrac{1-\gamma+\mu}{1-\mu}\Big]||T_{3}z-x_{m}||+\alpha_{m}\Big[\dfrac{1-\gamma+\mu}{1-\mu}\Big]||T_{3}z-x_{m}||+\alpha^{2}_{m}||z-x_{m}||\\
		&&+||(1-\alpha_{m}+\alpha^{2}_{m})T_{1}x_{m}+\alpha_{m}T_{2}x_{m}-\alpha^{2}_{m}T_{3}x_{m}-x_{m}||\\
		&\leq&(1-\alpha_{m}-\alpha^{2}_{m})||T_{3}z-x_{m}||+\alpha_{m}||T_{3}z-x_{m}||+\alpha^{2}_{m}||z-x_{m}||\\
		&&+||(1-\alpha_{m}+\alpha^{2}_{m})T_{1}x_{m}+\alpha_{m}T_{2}x_{m}-\alpha^{2}_{m}T_{3}x_{m}-x_{m}||~~(since~\dfrac{1-\gamma+\mu}{1-\mu}\leq1)\\
		&=&(1-\alpha^{2}_{m})||T_{3}z-x_{m}||+\alpha^{2}_{m}||z-x_{m}||\\
		&&+||(1-\alpha_{m}+\alpha^{2}_{m})T_{1}x_{m}+\alpha_{m}T_{2}x_{m}-\alpha^{2}_{m}T_{3}x_{m}-x_{m}||\\
		&<&(1-\alpha^{2}_{m})||T_{3}z-x_{m}||+\alpha^{2}_{m}\epsilon+\alpha^{2}_{m}\epsilon\\&=&(1-\alpha^{2}_{m})||T_{3}z-x_{m}||+2\alpha^{2}_{m}\epsilon
		\\&<&(1-\alpha^{2}_{m})||T_{3}z-x_{m}||+\alpha^{2}_{m}||T_{3}z-x_{m}||\\
		&=&||T_{3}z-x_{m}||,
	\end{eqnarray*}
	which is a contradiction. Hence, $z$ is a common fixed point of $T_{1}, T_{2}$ and $T_{3}~~~~~~~~~~~~~~~~~~~~~~~~~~~~~~~~~~~~\blacksquare$.
	\begin{lm} \label{lemma6}
			Let $C$ be a nonempty closed convex  subset of a  Banach space $X$. Let $T_{1}, T_{2}$ and $T_{3}$ be commuting  mappings   satisfying the condition $B_{\gamma,\mu}$ on $C$.  Let $\{\alpha_{n}\}$ be a sequence in $[0,\dfrac{1}{2}]$ converging to $0$. 
	Define a sequence $\{x_{n}\}$ in $C$ by $x_{0}\in C$ and
		\begin{equation}\label{equation10}
		x_{n+1}=\lambda(1-\alpha_{n}-\alpha^{2}_{n})T_{1}x_{n}+\lambda\alpha_{n}T_{2}x_{n}+\lambda\alpha^{2}_{n}T_{3}x_{n}+(1-\lambda)x_{n},
	\end{equation}
	for $\lambda\in[\gamma,1)-\{0\},~n=0,1,2,\cdots$.
		 Then
		 	\begin{equation}\label{eqn13}
		 \lim_{n\rightarrow\infty}||(1-\alpha_{n})T_{1}x_{n}+\alpha_{n}T_{2}x_{n}+\alpha^{2}_{n}T_{3}x_{n}-x_{n}|| =0.
		 \end{equation}
		  	
\end{lm}
\textbf{Proof.} Since  a sequence $\{x_{n}\}$ in $C$ defined as (\ref{equation10}) and assume that $z$ is a common fixed point of $T_{1}, T_{2}$ and $T_{2}  $,
then,   by  condition $B_{\gamma,\mu}$ we get
\begin{eqnarray*}
	\gamma||(1-\alpha_{n})T_{1}x_{n}+\alpha_{n}T_{2}x_{n}+\alpha_{n}^{2}T_{3}x_{n}-x_{n}||&\leq&\lambda||(1-\alpha_{n})T_{1}x_{n}+\alpha_{n}T_{2}x_{n}+\alpha_{n}^{2}T_{3}x_{n}-x_{n}||\\&=&||x_{n+1}-x_{n}||\\&\leq&||x_{n+1}-x_{n}||+\mu||x_{n}-T_{1}x_{n}||.	
\end{eqnarray*}
Implies that
\begin{eqnarray*}
	&&||(1-\alpha_{n+1}-\alpha_{n+1}^{2})T_{1}x_{n+1}+\alpha_{n+1}T_{2}x_{n+1}+\alpha_{n+1}^{2}T_{3}x_{n+1}-T_{1}x_{n}||
	\\&\leq&(1-\alpha_{n+1}-\alpha_{n+1}^{2})||T_{1}x_{n+1}-T_{1}x_{n}||+\alpha_{n+1}||T_{2}x_{n+1}-T_{1}x_{n}||
	+\alpha_{n+1}^{2}||T_{3}x_{n+1}-T_{1}x_{n}||\\&\leq&(1-\alpha_{n+1}-\alpha_{n+1}^{2})\Big[||T_{1}x_{n+1}-T_{1}z||+||T_{1}z-T_{1}x_{n}||\Big]+
	\alpha_{n+1}\Big[||T_{2}x_{n+1}-T_{2}z||+||T_{2}z-T_{1}x_{n}||\Big]
	\\&&+~~\alpha_{n+1}^{2}\Big[||T_{3}x_{n+1}-T_{3}z||+||T_{3}z-T_{1}x_{n}||\Big]
	\\&\leq&
	(1-\alpha_{n+1}-\alpha_{n+1}^{2})\Big[(1-\gamma)||x_{n+1}-z||+\mu(||z-T_{1}x_{n+1}||+||x_{n+1}-T_{1}z||)+(1-\gamma)||x_{n}-z||\\&&+\mu(||z-T_{1}x_{n}||+||x_{n}-T_{1}z||)\Big]
	+\alpha_{n+1}\Big[(1-\gamma)||x_{n+1}-z||+\mu(||z-T_{2}x_{n+1}||+||x_{n+1}-T_{2}z||)+\\&&
	(1-\gamma)||x_{n}-z||+\mu(||z-T_{1}x_{n}||+||x_{n}-T_{2}z||)\Big]
	+\alpha_{n+1}^{2}\Big[(1-\gamma)||x_{n+1}-z||+\mu(||z-T_{3}x_{n+1}||+||x_{n+1}-T_{3}z||)+\\&&
	(1-\gamma)||x_{n}-z||+\mu(||z-T_{1}x_{n}||+||x_{n}-T_{3}z||)\Big]
	\\&\leq&
	(1-\alpha_{n+1}-\alpha_{n+1}^{2})\Big[(1-\gamma)||x_{n+1}-z||+\mu(||z-x_{n+1}||+||x_{n+1}-z||)+(1-\gamma)||x_{n}-z||\\&&+\mu(||z-x_{n}||+||x_{n}-z||)\Big]
	+\alpha_{n+1}\Big[(1-\gamma)||x_{n+1}-z||+\mu(||z-x_{n+1}||+||x_{n+1}-z||)+\\&&
	(1-\gamma)||x_{n}-z||+\mu(||z-x_{n}||+||x_{n}-z||)\Big]
	+\alpha_{n+1}^{2}\Big[(1-\gamma)||x_{n+1}-z||+\mu(||z-x_{n+1}||+||x_{n+1}-z||)+\\&&
	(1-\gamma)||x_{n}-z||+\mu(||z-x_{n}||+||x_{n}-z||)\Big]
		\\&\leq&
	(1-\alpha_{n+1}-\alpha_{n+1}^{2})\Big[(1-\gamma)||x_{n+1}-z||+\gamma||z-x_{n+1}||+(1-\gamma)||x_{n}-z||+\gamma||z-x_{n}||\Big]
	+\\&&\alpha_{n+1}\Big[(1-\gamma)||x_{n+1}-z||+\gamma||z-x_{n+1}||+
	(1-\gamma)||x_{n}-z||+\gamma||z-x_{n}||\Big]
	+\\&&\alpha_{n+1}^{2}\Big[(1-\gamma)||x_{n+1}-z||+\gamma||z-x_{n+1}||+
	(1-\gamma)||x_{n}-z||+\gamma||z-x_{n}||\Big]
	\\&=&
	(1-\alpha_{n+1}-\alpha_{n+1}^{2})\Big[||x_{n+1}-z||+||x_{n}-z||\Big]
	+\alpha_{n+1}\Big[||x_{n+1}-z||+
	||x_{n}-z||\Big]
	+\alpha_{n+1}^{2}\Big[||x_{n+1}-z||+
	||x_{n}-z||\Big]\\&=&
	||x_{n+1}-z||+
	||x_{n}-z||	
\end{eqnarray*}
Now from the following we have that, 
\[||(1-\alpha_{n}-\alpha_{n}^{2})T_{1}x_{n}+\alpha_{n}T_{2}x_{n}+\alpha_{n}^{2}T_{3}x_{n}-x_{n}||\leq(1-\alpha_{n}-\alpha_{n}^{2})||T_{1}x_{n}-x_{n}||+\alpha_{n}||T_{2}x_{n}-x_{n}||+\alpha_{n}^{2}||T_{3}x_{n}-x_{n}||.\]
Taking $n\rightarrow\infty$ on both side of the above inequality and using Proposition \ref{proposition2} we get
\[\lim_{n\rightarrow\infty}||(1-\alpha_{n}-\alpha_{n}^{2})T_{1}x_{n}+\alpha_{n}T_{2}x_{n}+\alpha_{n}^{2}T_{3}x_{n}-x_{n}||=0.\]
		\begin{theorem}\label{theorem2}
		Let $C$ be a compact and convex subset of a  Banach space $X$. Let $T_{1}, T_{2}$ and $T_{3}$ be commuting  mappings   satisfying the condition $B_{\gamma,\mu}$ on $C$. Let $\{\alpha_{n}\}$ be a sequence in $[0,\dfrac{1}{2}]$ converging to $0$. 
		For $x_{0}\in C$, let $\{x_{n}\}$ be a sequence in $C$ as defined by (\ref{equation10}).
		 Then $\{ x_{n}  \}$ converges strongly to a common fixed point $T_{1}, T_{2}$ and $T_{3}$.	
	\end{theorem}
\textbf{Proof.} 
Since $C$ is compact, there exists a subsequence
$\{x_{n_{j}}\}$ of $\{x_{n}\}$ and $z\in C$ such that $\{x_{n_{j}}\}$ converges to $z\in C$ By (Lemma \ref{lemma4})\\
Now by Proposition \ref{proposition1} (ii), for $\gamma=\frac{\theta}{2},~\theta\in[0,1]$,
\begin{eqnarray*}
	&&\gamma||x_{n_{j}}-T_{1}x_{n_{j}}||\leq||x_{n_{j}}-z||\\
	&&\Rightarrow\gamma||x_{n_{j}}-T_{1}x_{n_{j}}||\leq||x_{n_{j}}-z||+\mu||z-T_{1}z||.
\end{eqnarray*}
So by condition $B_{\gamma,\mu}$ and by using (\ref{eqn7}) we have
\begin{eqnarray*}
	||x_{n_{j}}-T_{1}z||&\leq&||x_{n_{j}}-T_{1}x_{n_{j}}||+||T_{1}x_{n_{j}}-T_{1}z||\\&\leq&||x_{n_{j}}-T_{1}x_{n_{j}}||+(1-\gamma)||x_{n_{j}}-z||+\mu(||x_{n_{j}}-T_{1}z||+||z-T_{1}x_{n_{j}}||)
	\\&\leq&||x_{n_{j}}-T_{1}x_{n_{j}}||+(1-\gamma)||x_{n_{j}}-z||+\mu(||x_{n_{j}}-T_{1}z||+||z-x_{n_{j}}||+||x_{n_{j}}-T_{1}x_{n_{j}}||).	
\end{eqnarray*}
So taking $n_{j}\rightarrow\infty$ and using Proposition \ref{proposition2}, we get
\begin{eqnarray*}
	&&(1-\mu)||z-T_{1}z||\leq 0\\
	&&\Rightarrow Tz=z,	
\end{eqnarray*}
since $\mu\neq1$, shows that $z$ is a fixed point $T_{1}$. Therefore, by applying Lemma \ref{lemma5} and Theorem \ref{theorem1},  $z$ is a common fixed point of $T_{1}$ and $T_{2}$.\\
	Then by using (\ref{equation10}) and  (\ref{equation7}), we have the following:
	\begin{eqnarray*}
		||x_{n+1}-z||&=&||\lambda(1-\alpha_{n}-\alpha^{2}_{n})T_{1}x_{n}+\lambda\alpha_{n}T_{2}x_{n}+\lambda\alpha^{2}_{n}T_{3}x_{n}+(1-\lambda)x_{n}-z||\\&\leq&
		\lambda(1-\alpha_{n}-\alpha^{2}_{n})||T_{1}x_{n}-z||+\lambda\alpha_{n}||T_{2}x_{n}-z||+\lambda\alpha_{n}^{2}||T_{3}x_{n}-z||+(1-\lambda)||x_{n}-z||
		\\&\leq&
		\lambda(1-\alpha_{n}-\alpha^{2}_{n})\Big[(1-\gamma)||x_{n}-z||+\mu(|x_{n}-T_{1}z||+||z-T_{1}x_{n}||)\Big]+\\&&\lambda\alpha_{n}\Big[(1-\gamma)||x_{n}-z||+\mu(|x_{n}-T_{2}z||+||z-T_{2}x_{n}||)\Big]+\\&&\lambda\alpha_{n}^{2}\Big[(1-\gamma)||x_{n}-z||+\mu(|x_{n}-T_{3}z||+||z-T_{3}x_{n}||)\Big]+(1-\lambda)||x_{n}-z||
		\\&\leq&
		\lambda(1-\alpha_{n}-\alpha^{2}_{n})\Big[(1-\gamma)||x_{n}-z||+\mu(|x_{n}-z||+||z-x_{n}||)\Big]+\\&&\lambda\alpha_{n}\Big[(1-\gamma)||x_{n}-z||+\mu(|x_{n}-z||+||z-x_{n}||)\Big]+\\&&\lambda\alpha_{n}^{2}\Big[(1-\gamma)||x_{n}-z||+\mu(|x_{n}-z||+||z-x_{n}||)\Big]+(1-\lambda)||x_{n}-z||
		\\&\leq&
		\lambda(1-\alpha_{n}-\alpha^{2}_{n})||x_{n}-z||+\lambda\alpha_{n}||x_{n}-z||+\lambda\alpha_{n}^{2}||x_{n}-z||+(1-\lambda)||x_{n}-z||~~~(sine~2\mu\leq\gamma)\\&=&||x_{n}-z||. 
	\end{eqnarray*}
	Thus, $\{||x_{n}-z||\}$ is a monotonically decreasing non-negative sequence of real numbers, and will converges to some real numbers. 
	Therefore, by similar arguments of  Lemma \ref{lemma4} 
	the sequence $\{x_{n}\}$ converges strongly to the common fixed point of $T_{1},T_{2}$ and $T_{3}.~~~~~~~~~~~~~~~~~~~~~~~~~~~~~~~\blacksquare$
		\begin{lm}\label{lemma7}
			 Let $C$ be a nonempty closed convex  subset of a  Banach space $X$ Banach space $X$. Let $l\in\mathbb{N}$ with $l\geq2$ and let $T_{1},T_{2},\cdots, T_{l}$  be commuting mappings  satisfying the condition $B_{\gamma,\mu}$ on $C$. Let $\{\alpha_{n}\}$ be a sequence in $(0,\dfrac{1}{2})$ converging to $0$, and let $\{x_{n}\}$ be a sequence in $C$ such that $\{x_{n}\}$  converges strongly to some $z\in C$.
		Then $z$   a common fixed point of $T_{1},T_{2}\cdots T_{l}$.  
	\end{lm}
	\textbf{Proof.} We will prove this lemma by induction. We have already proved the conclusion in the case $l=2,3$. Fix $l\in\mathbb{N},~l\geq4$. We assume the conclusion holds for every integer less than $m$ and grater than $1$.
	
	Let $T_{1},T_{2},\cdots, T_{l}$  be commuting mappings
	satisfying the condition $B_{\gamma,\mu}$ on $C$. Then by
	Lemma \ref{lemma5} and the induction assumption, we have that $z$ is a common fixed point of
	$T_{1},T_{2},\cdots, T_{l-1}$.
	And we note that;
	\begin{equation}\label{equation12}
		T_{k}\circ T_{l}z=T_{l}\circ T_{k}z=T_{l}z,~with~1\leq k< l,~\forall k\in\mathbb{N}.
	\end{equation}
	We assume that $z$ is not a fixed point of $T_{m}$.  Put \[\epsilon=\dfrac{||T_{l}z-z||}{3}>0.\]
	Then there exists $m\in\mathbb{N}$ such that 
	\[||x_{m}-z||<\epsilon,~\frac{||(1-\sum_{k=1}^{l-1}\alpha^{k}_{m})T_{1}x_{m}+\sum_{k=2}^{l} \alpha_{m}^{k-1}T_{k}x_{m}-x_{m}||}{\alpha^{l-1}_{m}}< \epsilon.\]
	Since 
	\begin{eqnarray*}
		3\epsilon=||T_{m}z-z|| &\leq&||T_{l}z-x_{m}||+||x_{m}-z||\\
		&<&||T_{l}z-x_{m}||+\epsilon\\ \Rightarrow 2\epsilon &<&||T_{l}z-x_{m}||.
	\end{eqnarray*}
	So, we obtain the following
	\begin{eqnarray*}
		||T_{l}z-x_{m}||&\leq & ||T_{l}z-(1-\sum_{k=1}^{l-1}\alpha_{m}^{k})T_{1}x_{m}-\sum_{k=2}^{l}\alpha_{m}^{k-1}T_{k}x_{m}||\\
		&&+||(1-\sum_{k=1}^{l-1}\alpha_{m}^{k})T_{1}x_{m}+\sum_{k=2}^{l}\alpha_{m}^{k-1}T_{k}x_{m}-x_{m}||\\&\leq&
		(1-\sum_{k=1}^{l-1}\alpha_{m}^{k})||T_{l}z-T_{1}x_{m}||+\sum_{k=2}^{l-1}\alpha_{m}^{k-1}||T_{l}z-T_{k}x_{m}||+\alpha_{m}^{l-1}||T_{l}z-T_{l}x_{m}||\\
		&&+||(1-\sum_{k=1}^{l-1}\alpha_{m}^{k})T_{1}x_{m}+\sum_{k=2}^{l}\alpha_{m}^{k-1}T_{k}x_{m}-x_{m}||
		\\&=&
		(1-\sum_{k=1}^{l-1}\alpha_{m}^{k})||T_{l}\circ T_{k} z-T_{1}x_{m}||+\sum_{k=2}^{l-1}\alpha_{m}^{k-1}||T_{l}\circ T_{k} z-T_{k}x_{m}||+\alpha_{m}^{l-1}||T_{l}z-T_{l}x_{m}||\\
		&&+||(1-\sum_{k=1}^{l-1}\alpha_{m}^{k})T_{1}x_{m}+\sum_{k=2}^{l}\alpha_{m}^{k-1}T_{k}x_{m}-x_{m}||
		\\&\leq&
		(1-\sum_{k=1}^{l-1}\alpha_{m}^{k})\Big((1-\gamma)||T_{l}z-x_{m}||+\mu(||T_{l} z-x_{m}+||x_{m}-T_{l}z||)\Big)+\\&&\sum_{k=2}^{l-1}\alpha_{m}^{k-1}\Big((1-\gamma)||T_{l}z-x_{m}||+\mu(||T_{k} z-x_{m}+||x_{m}-T_{l}z||)\Big)+\alpha_{m}^{l-1}||z-x_{m}||\\
		&&+||(1-\sum_{k=1}^{l-1}\alpha_{m}^{k})T_{1}x_{m}+\sum_{k=2}^{l}\alpha_{m}^{k-1}T_{k}x_{m}-x_{m}||
			\\&\leq&
		(1-\sum_{k=1}^{l-1}\alpha_{m}^{k})\Big[\dfrac{(1-\gamma+\mu)}{1-\mu}\Big]||T_{l}z-x_{m}||+\sum_{k=2}^{l-1}\alpha_{m}^{k-1}\Big[\dfrac{(1-\gamma+\mu)}{1-\mu}\Big]||T_{l}z-x_{m}||+\alpha_{m}^{l-1}||z-x_{m}||\\
		&&+||(1-\sum_{k=1}^{l-1}\alpha_{m}^{k})T_{1}x_{m}+\sum_{k=2}^{l}\alpha_{m}^{k-1}T_{k}x_{m}-x_{m}||
			\\&\leq&
		(1-\sum_{k=1}^{l-1}\alpha_{m}^{k})||T_{l}z-x_{m}||+\sum_{k=2}^{l-1}\alpha_{m}^{k-1}||T_{l}z-x_{m}||+\alpha_{m}^{l-1}||z-x_{m}||\\
		&&+||(1-\sum_{k=1}^{l-1}\alpha_{m}^{k})T_{1}x_{m}+\sum_{k=2}^{l}\alpha_{m}^{k-1}T_{k}x_{m}-x_{m}||
			\\&=&
		||T_{l}z-x_{m}||+\alpha_{m}^{l-1}||z-x_{m}||+||(1-\sum_{k=1}^{l-1}\alpha_{m}^{k})T_{1}x_{m}+\sum_{k=2}^{l}\alpha_{m}^{k-1}T_{k}x_{m}-x_{m}||
		\\&<&
		||T_{l}z-x_{m}||+\alpha_{m}^{l-1}\epsilon+\alpha_{m}^{l-1}\epsilon
		\\&=&
		(1-\alpha_{m}^{l-1})||T_{l}z-x_{m}||+2\alpha_{m}^{l-1}\epsilon
		\\&<&
		(1-\alpha_{m}^{l-1})||T_{l}z-x_{m}||+\alpha_{m}^{l-1}||T_{l}z-x_{m}||
		\\&=&
		||T_{l}z-x_{m}||,
	\end{eqnarray*}
	which is a contradiction. Hence, $z$ is a common fixed point of $T_{1}, T_{2},\cdots,T_{l},~l\in\mathbb{N}$. $~~~~~~~~~~~~~~~~~~~~~~~~~~~~~~~~~~~~~~~\blacksquare$.
	\begin{lm} \label{lemma8}
		Let $C$ be a nonempty closed convex and bounded subset of a  Banach space $X$. Let $T_{1}, T_{2},\cdots,T_{l},~l\in\mathbb{N}$ be commuting  mappings   satisfying the condition $B_{\gamma,\mu}$ on $C$.  Let $\{\alpha_{n}\}$ be a sequence in $[0,\dfrac{1}{2}]$ converging to $0$. 
		Define a sequence $\{x_{n}\}$ in $C$ by $x_{0}\in C$ and
		\begin{equation}\label{equation13}
		x_{n+1}=\lambda\Big(1-\sum_{k=1}^{l-1}\alpha_{n}^{k}\Big)T_{1}x_{n}+\lambda\sum_{k=2}^{l}\alpha_{n}^{k-1}T_{k}x_{n}+(1-\lambda)x_{n},
		\end{equation}
		for $\lambda\in[\gamma,1)-\{0\},~n=0,1,2,\cdots$.
		Then
		\begin{equation}\label{eqn17}
			\lim_{n\rightarrow\infty}||\Big(1-\sum_{k=1}^{l-1}\alpha_{n}^{k}\Big)T_{1}x_{n}+\sum_{k=2}^{l}\alpha_{n}^{k-1}T_{k}x_{n}-x_{n}|| =0.
		\end{equation}
		 	
	\end{lm}
	\textbf{Proof.} Since   and a sequence $\{x_{n}\}$ in $C$ defined as (\ref{equation13}) and assume that $z$ is a common fixed point of $T_{1}, T_{2},\cdots T_{l}$ 	then,   by  condition $B_{\gamma,\mu}$ we get
	\begin{eqnarray*}
	\gamma||\Big(1-\sum_{k=1}^{l-1}\alpha_{n}^{k}\Big)T_{1}x_{n}+\sum_{k=2}^{l}\alpha_{n}^{k-1}T_{k}x_{n}-x_{n}||&\leq&\lambda||\Big(1-\sum_{k=1}^{l-1}\alpha_{n}^{k}\Big)T_{1}x_{n}+\sum_{k=2}^{l}\alpha_{n}^{k-1}T_{k}x_{n}-x_{n}||\\
	&=&||x_{n+1}-x_{n}||\\&\leq&|x_{n+1}-x_{n}||+\mu||x_{n}-T_{1}x_{n}||.	
	\end{eqnarray*}
It implies that
	\begin{eqnarray*}
&&	||\Big(1-\sum_{k=1}^{l-1}\alpha_{n+1}^{k}\Big)T_{1}x_{n+1}+\sum_{k=2}^{l-1}\alpha_{n+1}^{k-1}T_{k}x_{n+1}-T_{1}x_{n}||\\&\leq&
\Big(1-\sum_{k=1}^{l-1}\alpha_{n+1}^{k}\Big)||T_{1}x_{n+1}-T_{1}x_{n}||+	\sum_{k=2}^{l}\alpha_{n+1}^{k-1}||T_{k}x_{n+1}-T_{1}x_{n}||+
\\&\leq&
\Big(1-\sum_{k=1}^{l-1}\alpha_{n+1}^{k}\Big)\Big[||T_{1}x_{n+1}-T_{1}z||+||T_{1}-T_{1}x_{n}||\Big]+	\sum_{k=2}^{l}\alpha_{n+1}^{k-1}\Big[||T_{k}x_{n+1}-T_{k}z||+||T_{k}z-T_{1}x_{n}||\Big]
\\&\leq&.
\Big(1-\sum_{k=1}^{l-1}\alpha_{n+1}^{k}\Big)\Big[(1-\gamma)||x_{n+1}-z||+\mu(||z-T_{1}x_{n+1}||+||x_{n+1}-T_{1}z||)+(1-\gamma)||x_{n}-z||+\\&&\mu(||z-T_{1}x_{n}||+||x_{n}-T_{1}z||)\Big]
+\sum_{k=2}^{l}\alpha_{n+1}^{k-1}\Big[(1-\gamma)||x_{n+1}-z||+\mu(||z-T_{k}x_{n+1}||+\\&&||x_{n+1}-T_{k}z||)+(1-\gamma)||x_{n}-z||+\mu(||z-T_{1}x_{n}||+||x_{n}-T_{k}z||)\Big]
\\&\leq&.
\Big(1-\sum_{k=1}^{l-1}\alpha_{n+1}^{k}\Big)\Big[(1-\gamma)||x_{n+1}-z||+\mu(||z-x_{n+1}||+||x_{n+1}-z||)+(1-\gamma)||x_{n}-z||+\\&&\mu(||z-x_{n}||+||x_{n}-z||)\Big]
+\sum_{k=2}^{l}\alpha_{n+1}^{k-1}\Big[(1-\gamma)||x_{n+1}-z||+\mu(||z-x_{n+1}||+\\&&||x_{n+1}-z||)+(1-\gamma)||x_{n}-z||+\mu(||z-x_{n}||+||x_{n}-z||)\Big]
\\&\leq&.
\Big(1-\sum_{k=1}^{l-1}\alpha_{n+1}^{k}\Big)\Big[(1-\gamma)||x_{n+1}-z||+\gamma||x_{n+1}-z||+(1-\gamma)||x_{n}-z||+\gamma||x_{n}-z||\Big]
+\\&&\sum_{k=2}^{l}\alpha_{n+1}^{k-1}\Big[(1-\gamma)||x_{n+1}-z||+\gamma||x_{n+1}-z||+(1-\gamma)||x_{n}-z||+\gamma||x_{n}-z||\Big]
\\&=&
\Big(1-\sum_{k=1}^{l-1}\alpha_{n+1}^{k}\Big)\Big[||x_{n+1}-z||+||x_{n}-z||\Big]
+\sum_{k=2}^{l}\alpha_{n+1}^{k-1}\Big[||x_{n+1}-z||+||x_{n}-z||\Big]\\
&=&||x_{n+1}-z||+||x_{n}-z||
.
\end{eqnarray*}
Then, we have 
\[||\Big(1-\sum_{k=1}^{l-1}\alpha_{n}^{k}\Big)T_{1}x_{n}+\sum_{k=2}^{l}\alpha_{n}^{k-1}T_{k}x_{n}-x_{n}||\leq \Big(1-\sum_{k=1}^{l-1}\alpha_{n}^{k}\Big)||T_{1}x_{n}-x_{n}||+\sum_{k=2}^{l}\alpha_{n}^{k-1}||T_{k}x_{n}-x_{n}||\]
Now by taking $n\rightarrow\infty$ and applying Proposition \ref{proposition2}, we get
\[\lim_{n\rightarrow\infty}||\Big(1-\sum_{k=1}^{l-1}\alpha_{n}^{k}\Big)T_{1}x_{n}+\sum_{k=2}^{l}\alpha_{n}^{k-1}T_{k}x_{n}-x_{n}|| =0.\]
\begin{theorem}\label{theorem3}
	Let $C$ be a compact and convex subset of a  Banach space $X$. Let $T_{1}, T_{2},\cdots,T_{l}$ be commuting  mappings   satisfying the condition $B_{\gamma,\mu}$ on $C$. Let $\{\alpha_{n}\}$ be a sequence in $[0,\frac{1}{2}]$ converging to $0$. 
	For $x_{0}\in C$, let $\{x_{n}\}$ be a sequence in $C$ as defined by (\ref{equation13}).
	Then $\{ x_{n}  \}$ converges strongly to a common fixed point of $T_{1}, T_{2},\cdots,T_{l}$.	
\end{theorem}
\textbf{Proof.}  By using Theorem \ref{theorem1} and induction assumption, assume that $z$ is a common fixed point of $T_{1},T_{2},\cdots,T_{l}$. Then by using (\ref{equation13}) and (\ref{equation7}), we have
\begin{eqnarray*}
||x_{n+1}-z||&=&||\lambda\Big(1-\sum_{k=1}^{l-1}\alpha_{n}^{k}\Big)T_{1}x_{n}+\lambda\sum_{k=2}^{l}\alpha_{n}^{k-1}T_{k}x_{n}+(1-\lambda)x_{n}-z||\\&\leq&
\lambda\Big(1-\sum_{k=1}^{l-1}\alpha_{n}^{k}\Big)||T_{1}x_{n}-z||+\lambda\sum_{k=2}^{l}\alpha_{n}^{k-1}||T_{k}x_{n}-z||+(1-\lambda)||x_{n}-z||
\\&\leq&\lambda\Big(1-\sum_{k=1}^{l-1}\alpha_{n}^{k}\Big)\Big[(1-\gamma)||x_{n}-z||+\mu(||x_{n}-T_{1}z||+||z-T_{1}x_{n}||)\Big]+\\
&&\lambda\sum_{k=2}^{l}\alpha_{n}^{k-1}
\Big[(1-\gamma)||x_{n}-z||+\mu(||x_{n}-T_{k}z||+||z-T_{k}x_{n}||)\Big]+(1-\gamma)||x_{n}-z||
\\&\leq&\lambda\Big(1-\sum_{k=1}^{l-1}\alpha_{n}^{k}\Big)\Big[\dfrac{1-\gamma+\mu}{1-\mu}\Big]||x_{n}-z||+\lambda\sum_{k=2}^{l}\alpha_{n}^{k-1}
\Big[\dfrac{1-\gamma+\mu}{1-\mu}\Big]||x_{n}-z||+(1-\gamma)||x_{n}-z||
\\&\leq&\lambda\Big(1-\sum_{k=1}^{l-1}\alpha_{n}^{k}\Big)||x_{n}-z||+\lambda\sum_{k=2}^{l}\alpha_{n}^{k-1}
||x_{n}-z||+(1-\gamma)||x_{n}-z||
\\&=&\lambda||x_{n}-z||+(1-\gamma)||x_{n}-z||
\\&=&||x_{n}-z||.
\end{eqnarray*} 
Thus, $\{||x_{n}-z||\}$ is a monotonically decrease non-negative sequence of real numbers, and will converge to some real number. Therefore, by Theorem \ref{theorem2} and induction assumption, the sequence $\{x_{n}\}$ converges strongly to the common fixed point of $T_{1},T_{2},\cdots,T_{l}.~~~~~~~~~~~~~~~~~~~~~~~~~~~~~~~~~~~~~~~~~~~\blacksquare$ 
	\section{Conclusion}
In this paper, our purpose was to extend the classic approach of fixed point searching from \cite{suzuki2005strong}, by taking the condition $B_{\gamma,\mu}$ and the process by properly using the finite number of operators in to its structure. In other words, we brought together both the Suziki's strong nonexpansive mappings and the  Patir et al. \cite{patir2018some} condition $B_{\gamma,\mu}$ under the same iteration process. Under the resulted iteration process, we proved the existence of a common fixed point for finite number of commutative mapping satisfies the condition $B_{\gamma,\mu}$. We proves several strong convergence results and existence of the common fixed point for the sequence of approximation by those iteration.
	\bibliographystyle{plain}
	\bibliography{geze2}
\end{document}